%% file: main.tex
\documentclass[12pt]{article}

\usepackage[english]{babel}
\usepackage[utf8x]{inputenc}
\usepackage{pdfpages}
\usepackage{amsfonts,dsfont}
\usepackage{amsmath}
\usepackage{amssymb}
\usepackage{amsthm} 
\usepackage{bm}
\usepackage{hyperref}
\usepackage[linesnumbered]{algorithm2e}
\usepackage{enumitem}
\usepackage[capitalise]{cleveref}

\usepackage{todonotes}
\usepackage[all]{xy} 

\makeatletter \let\cl@part\relax \makeatother

\usepackage[top=1.5cm, bottom=1.5cm, left=1.5cm, right=1.5cm]{geometry}
\usepackage{xcolor}
\usepackage{tikz}
\usetikzlibrary{calc}
\usetikzlibrary{intersections}
\tikzset{offset/.style={to path={%
    -- ($(\tikztostart)!#1cm!(\tikztotarget)$)}},
         offset/.default=1}
\tikzset{>=latex}

\usepackage{subcaption}
\usepackage{tikz-3dplot}

\Crefname{ass}{Assumption}{Assumptions}

\newtheorem{theorem}{Theorem}
\newtheorem{lemma}[theorem]{Lemma}
\newtheorem{remark}[theorem]{Remark}

\newtheorem{defi}[theorem]{Definition}

\theoremstyle{remark}

\newcommand{\GGAPM}{G²APM}
\newcommand{\GGAPMspace}{G²APM }

\title{Generalized adaptive partition-based method for
two-stage stochastic linear programs : convergence and generalisation}
\author{Maël Forcier, Vincent Leclère}
\input{commandes}
\begin{document}

\maketitle

\abstract{
	Adaptive Partition-based Methods (APM) are numerical methods to solve two-stage stochastic linear problems (2SLP). 
The core idea is to iteratively construct an \emph{adapted} partition of the space of alea in order 
to aggregate scenarios while conserving the true value of the cost-to-go for the current first-stage control.
Relying on the normal fan of the dual admissible set, we extend the classical and generalized APM method by 
i) extending the method to almost arbitrary 2SLP, 
ii) giving a necessary and sufficient condition for a partition to be adapted even for non-finite distribution, and 
iii) proving the convergence of the method.
We give some additional insights by linking APM to the L-shaped algorithm.
}

\section{Introduction}

Stochastic programming is a powerful modeling paradigm for optimization under uncertainty that has found many applications in energy, logistics or finance (see e.g.\ \cite{wallace2005applications}).
Two stage linear stochastic programs (2SLP) constitute an important class of stochastic programs.
They have been thoroughly studied, see e.g.\ \cite{
birge2011introduction
} and references therein.
One reason for this interest is the availability of efficient linear solvers and the use of dedicated algorithms leveraging the special structure of linear stochastic programs (\cite{van1969shaped,birge1985decomposition}). 

\subsection{Setting}
Let $(\Omega,\cA,\PP)$ be a probability space and $\va \xi=(\va T,\va h)\in L_1(\Omega,\cA,\PP;\Xi)$ an integrable random variable with values in $\Xi:=\RR^{l\times n} \times \RR^l $.
 We consider the following problem 2-stage stochastic linear problem with fixed recourse (\ie\ $W$ is deterministic): 

\begin{equation}
		\min_{x\in \RR_+^n} \; \Ba{c^\top x + \besp{Q(x,\va \xi)} \; | \; Ax=b} 
\tag{2SLP}
\label{eq:2SLP}
\end{equation}
where the recourse cost is
\begin{equation}
	Q(x,\xi):= \; \min_{y \in \RR_+^m}  \; 
	\Ba{q^\top y \; | \; 
		 T x + W y = h }
	\label{eq:2nd_stage}
\end{equation}
The dual formulation of the recourse cost is
\begin{equation}
	Q^D(x,\xi):= \; \max_{\dualvar \in \RR^l} \; \Ba{(h - T x)^\top \dualvar 
		\; | \; W^\top \dualvar \leq q }\, .
		\label{eq:dual_2nd_stage}
\end{equation}
	We define 
	\begin{equation}
		X:= \{x \in \RR_+^n\; |\; Ax=b\},
	\end{equation}
	\begin{equation}
		D := \{ \dualvar \in \RR^l \; |\;  W^\top \dualvar \leq q \}
		\label{eq:defi_dual_set}
	\end{equation}

In the rest of the paper, we assume that $D \neq \emptyset$
which implies strong duality: $Q(x,\xi)=Q^D(x,\xi)$.

\subsection{Contribution and literature review}

Most results for 2SLP with continuous distributions rely on discretizing the distributions.
The Sample Average Approximation (SAA) method  samples the costs and constraints.
It relies on probabilistic results based on a uniform law of large number to give statistical guarantees,
see \ \cite[Chap. 5]{shapiro2014lectures} for details.
Obtaining an approximation with satisfying guarantees requires a large number of scenarios.
Otherwise, when the support of the random variables are simplices, we can leverage convexity inequalities (like Jensen's and Edmundson-Madansky's) or moments inequalities
to construct finite scenario trees such that the discretized problem 
yields upper or lower bound of the continuous one (see \emph{e.g.} \cite{kuhn2006generalized,edirisinghe1994bounds}).

In each of these approaches, we solve 
an approximate version of the stochastic program, with or without  guarantees.
In any case, the number of scenarios increase the numerical burden of 2SLP. 

In order to alleviate the computations, we can use scenario reduction methods.
Some are based on heuristics, aiming at matching properties of the underlying distribution (e.g. matching moments), others are based on adequate distances on the scenario tree (see \cite{dupacova2003scenario, heitsch2007note}).
Recently Song and Luedtke \cite{song2015adaptive} proposed an exact scenario reduction method, 
which yields a reduced problem with the same value as the original.
This adaptive partition-based method (APM) is designed for finitely supported, fixed recourse, problems.
It consists in partitioning the space of alea, and
replacing each subset of scenarios by their weighted mean.
We say that a partition is \emph{adapted} to a first-stage control $x$
if the aggregated recourse problem has the same value as the recourse problem with the original distribution.
APM iteratively refines the partition so that, at each step, the new partition is
adapted to the current first-stage control.

	Van Ackooj, de Oliveira and Song \cite{van2018adaptive} improved the performance of APM by combining it with level decomposition methods ideas.
Finally, Ramirez-Pico and Moreno extended Song and Luedtke algorithm in \cite{ramirez2021generalized} to problems with continuous distribution for the right-hand side (and deterministic costs and matrices of constraints).
Their GAPM algorithm relies on a sufficient condition on the optimal dual variables for the partition to be adapted.
Nevertheless, finding a partition which satisfies this condition is left to the user, thus requiring an analytical approach for each particular problem.

Thanks to polyhedral geometry tools, we give a new necessary and sufficient condition for a partition to be adapted. 
It relies on the normal fan of the dual admissible set.
This allows us to give a general algorithm, called \GGAPM, that can be used for any 2SLP with fixed recourse, including problems with a general stochastic cost.
In addition, we get converging upper bounds  for free.

\GGAPMspace
 can also be extended to cases with general cost distribution, leaving only the recourse matrix $\va W$ with a finite support.

Finally, we show that \GGAPM, GAPM and APM can be interpreted as acceleration of L-shaped methods where we approximate the expected cost-to-go function by an exact tangent cone instead of a cutting hyperplane.
This allows us to prove the convergence of 
both GAPM and \GGAPMspace and provide an explicit bound on the number of steps.

\subsection{Structure of the paper}

\cref{sec:adapted_partition} present a generic framework for APM
and our \GGAPMspace algorithm in \cref{sec:adapted_partition}. 
By comparing it to known algorithms, we prove the convergence of \GGAPMspace in \cref{sec:comparison_convergence}.
Finally, \cref{sec:numerical} gives implementation details and present numerical results.

\section{Algorithm}
\label{sec:algo}
In this section we start by presenting a generic framework for APM 
algorithms, which depends on the choice of \emph{adapted partition}.
We proceed by giving a necessary and sufficient condition for a partition to be adapted, 
which leads to the definition of \GGAPM.

\subsection{Partition, refinements and APM framework}
\label{sub:partition}

We use $\besp{\cdot | P}$ to denote the conditional expectation $\besp{\cdot | \va \xi \in P}$ and $\bprob{P}$ to denote the probability $\bprob{\va \xi \in P}$.

We are considering partitions over the space of randomness, thus $\PP$-negligeable differences are irrelevant.
Consequently, we say that two sets of $E, F \subset\Xi$ are \emph{$\PP$-equivalent}, denoted $E \sim_\PP F$, 
if and only if they differ by a $\PP$-negligeable set 
\begin{equation}
	E \sim_\PP F\iff 
	\bprob{E\cap F}=\bprob{E}=\bprob{F},
\end{equation}
similarly we denote
\begin{equation}
	E \subset_\PP F \iff \bprob{E\cap F}=\bprob{E}.
\end{equation}
	Finally, we define a $\PP$-partition of $\Xi$ as the equivalence class of all partitions that are $\PP$-equivalent.

	Let $\cP$ and $\cR$ be two $\PP$-partitions of $\Xi$. We say that $\cP$ refines $\cR$, denoted $\cP \refleq_\PP \cR$, if
	\begin{equation}
		\forall P \in \cP, \; \exists R \in \cR, \quad P \subset_\PP R
	\end{equation}
The \emph{common refinement} of $\cP$ and $\cR$ is given by
	\begin{equation}
		\cP \wedge \cR =\{P \cap R \,|\, P \in \cP, R \in \cR\}.
	\end{equation}

\begin{defi}[Expected-cost-go of partition]
	Let $\cP$ be a (finite) $\PP$-partition of $\Xi \subset \RR^{l \times m}\times \RR^l$, then\footnote{We  use the convention that for every scalar or vector $z$, $z\times 0=0$. Thus, even if $\besp{\va \xi| P}$ might be undefined, it does not matter since the factor where it appears are multiplied by $\bprob{ P}=0$.}
	\begin{equation}
		\ectg_\cP(x):= \sum_{P \in \cP} \bprob{ P}Q\bp{x,\besp{\va \xi| P}} 
		\label{eq:defi_ectg_partition}
	\end{equation}

	We  say that a partition $\cP$ is adapted to $x$ if $\ectg_{\cP}(x)=\ectg(x)\defegal \besp{Q(x,\va \xi)}$.
\end{defi}

The following lemma show that, by convexity, a finer partition yields a higher expected cost-to-go function.
\begin{lemma}
	Let $\cP$ and $\cR$ two $\PP$-partitions of $\Omega$ then
	\begin{equation}
		\cP \refleq_\PP \cR \Rightarrow \ectg_{\cR} \leq \ectg_{\cP}
	\end{equation}
	\label{lem:monotony_ref}
	Moreover,
	\begin{equation}
		\label{eq:Monotonicity_of_V_P}
		\ectg_{\cP \wedge \cR} \geq\max(\ectg_\cP, \ectg_{\cR})
	\end{equation}
	In particular,
	\begin{equation}
		\label{eq:monotonicity_specified}
		Q\bp{\cdot,\besp{\va \xi}} \leq 
		\ectg_\cP
		\leq 
		\ectg .
	\end{equation}
\end{lemma}

\begin{proof}
By Jensen inequality and convexity of $\xi \mapsto Q(x,\xi)$, when $\bprob{R}>0$,
\begin{equation*}
	Q\bp{x,\besp{\va \xi|R}} \leq \sum_{P \in \cP} Q\bp{x,\besp{\va \xi|P\cap R}} \frac{\bprob{P \cap R}}{\bprob{R}}
\end{equation*}
Then, if $\cP\refleq_\PP \cR$,
	\begin{subequations}
		\begin{align*}
			\ectg_\cR(x) &=\sum_{R \in \cR} Q\bp{x,\besp{\va \xi|R}} \bprob{R} \\
			&\leq \sum_{P \in \cP}\sum_{R \in \cR} Q\bp{x,\besp{\va \xi|P\cap R}}\bprob{P \cap R} \\
			&= \sum_{P \in \cP} Q\bp{x,\besp{\va \xi|P}}\bprob{P} =\ectg_\cP(x)
		\end{align*}
	\end{subequations}
	The last line follows from the fact, that for $P \in \cP$, with $\bprob{P>0}$ and $\cP \refleq_\PP \cR$, 
	there exists a unique $R \in \cR$ such that $\bprob{P\cap R}=\bprob{P}$, all other $R \in \cR$ being such that $\bprob{P \cap R}=0$.

		\cref{eq:Monotonicity_of_V_P} is a direct consequence of $\ectg_{\cP \wedge \cP'} \geq \ectg_\cP$ and  $\ectg_{\cP \wedge \cP'} \geq \ectg_{\cP'}$.
	Thus, $\ectg_{\cP \wedge \cP'} \geq\max(\ectg_\cP, \ectg_{\cP'})$.
	Coupling this result with $\cP \refleq_\PP \{\Xi\}$ yield the left inequality of \cref{eq:monotonicity_specified} while the other can be found in \cite[Prop. 1]{ramirez2021generalized}.
\end{proof}

With those definitions we present a generic framework for APM methods.

\begin{algorithm}
	$k\leftarrow 0$, $z^U_{0}\leftarrow +\infty$, $z^L_{0}\leftarrow -\infty$, $\cP^{0}\leftarrow \{\Xi\}$ \;
	\While{$z^U_{k}-z^L_{k} >\varepsilon$}
	{
		$k\leftarrow k+1$\;
		Solve $z^L_k \leftarrow \min_{x \in X} c^\top x + \ectg_{\cP^{k-1}}(x)$  
		and let $x_k$ be an optimal solution \label{li:solve_part_lp}\;
		Choose a partition $\cP_{x_k}$ adapted to $x_k$ 
		\label{li:adapted_partition}
		\;
		$\cP^{k}\leftarrow  \cP^{k-1}\wedge \cP_{x_k}$ \label{li:common_refinement}\;
		$z^U_{k}\leftarrow \min\Bp{z^U_{k-1}, c^\top x_{k} + \ectg_{\cP^k}(x_k)}$ \label{li:upper_bound}\;
	}
	\caption{Generic framework for APM}
	\label{algo:GAPM}
\end{algorithm}

The key step of \cref{algo:GAPM}, in line~\ref{li:adapted_partition}, not yet specified, is the construction of an adapted partition 
$\cP_{x_k}$.
In this paper we study a specific APM, called \GGAPM, where we choose the partition $\cR_{x_k}$ defined in \cref{thm:ref_equiv_exact}. 

\subsection{Coarsest adapted partition}
\label{sec:adapted_partition}

In this section, we give the definition of a particular partition $\cR_x$ and prove that it is, in a generic case, the coarsest partition adapted to $x$.

When the distributions have finite support, \cite{song2015adaptive} 
characterized the partitions adapted to $x$.
Building on this result, a sufficient condition for continuous distribution can be found in \cite[Prop. 2]{ramirez2021generalized}.
We now prove that, for any distribution, a partition is adapted to $x$ if and only if it refines the collection $\bar \cR_x$. 
Unfortunately, $\bar \cR_x$ is not necessarily a partition, thus we also provide
a partition $\cR_x \preceq \bar \cR_x$ (see \Cref{fig:Ra_and_Rbarx} for an illustration).

	Recall that $D=\{\lambda\in \RR^l \;|\; W^\top \lambda \leq q\}$ and let  $\cN(D)$ be the normal fan of $D$ \ie\ the (finite) collection of all normal cones of $D$ (see \cite{forcier2021exact} for an introduction on normal fans and their use in stochastic programming).
	We denote by $\collecmax{\cN(D)}$ the collection of the maximal elements of $\cN(D)$ (full dimensional cones up-to lineality spaces).

\begin{theorem}

	Consider $x\in \RR^n$ and $N$ a cone of $\RR^m$. We define scenarios cells 
	$E_{N,x}$ and $\bar E_{N,x}$, subsets of $\Xi$, as
	\begin{subequations}
		\begin{align}
			E_{N,x}&:=\{ \xi \in \Xi \; |\; h-T x \in \relint(N)\} \\
			\bar E_{N,x}&:=\{ \xi \in \Xi \; |\; h-T x\in N\}
		\end{align}	
	\end{subequations}
	We define $\cR_x$ and $\bar \cR_x$ as
	\begin{subequations}
		\begin{align}
		\cR_x &:= \ba{ E_{N,x} \; |\; N \in \cN(D)} \\
		\bar \cR_x &:= \ba{ \bar E_{N,x} \; |\; N \in \collecmax{\cN(D)}} .
	\end{align}
	\end{subequations}
	Then,  
	\begin{subequations}
	\begin{align}
\cP \refleq_\PP \cR_x &\Longrightarrow 
\ectg_\cP(x)=\ectg(x) \\
\cP \refleq_\PP \bar \cR_x &\iff
\ectg_\cP(x)=\ectg(x). 
	\end{align}
	\label{eq:ref_equiv_exact}
\end{subequations}	
	\label{thm:ref_equiv_exact}
\end{theorem}
\begin{figure}[ht]
\def\scalecRx{0.45}

\def\raycone{6}
\def\opacitycone{0.2}
\def\opacitypoly{0.3}
\def\lengthcone{6.2}
\def\rayreccone{2}
\def\lengthreccone{1.9}
\def\colorcone{green!80!white}
\def\colorboundcone{green!70!blue}
\def\colorreccone{red}
\def\colorpoly{orange}
\def\colorboundpoly{orange}
\def\colorbigpoly{blue}
\def\colorcuttingplane{yellow!60!white}
\def\colorcomplex{red}
\def\colortriangulation{black}
\def\origin{0,0,0}
\def\colorvrtx{green!80!black}

\def\decal{0.5}

\begin{subfigure}{.2\textwidth}
    
\begin{tikzpicture}[scale=\scalecRx]

\pgfmathsetmacro\anga{0};
\pgfmathsetmacro\angb{45};
\pgfmathsetmacro\angc{90};
\pgfmathsetmacro\angd{150};
\pgfmathsetmacro\ange{170};

\pgfmathparse{\anga>\angb}%
        \ifnum\pgfmathresult=1
                \pgfmathsetmacro\anga{\anga+360}%
        \fi
\fill[\colorcone,opacity=\opacitycone] ($(\anga:\raycone)+(\decal,0)$) -- (\decal,0) --  ($(\angb:\raycone)+(\decal,0)$)  arc (\angb:\anga:\raycone);

\pgfmathparse{\angb>\angc}%
        \ifnum\pgfmathresult=1
                \pgfmathsetmacro\angb{\angb+360}%
        \fi
\fill[\colorcone,opacity=\opacitycone] ($(\angb:\raycone)+(0,\decal)$) -- (0,\decal) -- ($(\angc:\raycone)+(0,\decal)$)  arc (\angc:\angb:\raycone);

\draw[thick,->,\colorboundcone,name path=halflinea] (\decal, - \decal) -- ($(\anga:\lengthcone)+(\decal,-\decal)$);
\draw[thick,->,\colorboundcone,name path=halflineb] (\decal,\decal) -- ($(\angb:\lengthcone)+(\decal,\decal)$);
\draw[thick,->,\colorboundcone,name path=halflinec] (-\decal,\decal) -- ($(\angc:\lengthcone)+(-\decal,\decal)$);

\draw[\colorvrtx] (0,0)node{$\bullet$};

\draw[dotted,thin,\colorboundcone] ($(\anga:\raycone)+(\decal,0)$) -- (\decal,0);
\draw[dotted,thin,\colorboundcone] ($(\angb:\raycone)+(0,\decal)$) -- (0,\decal);
\draw[dotted,thin,\colorboundcone] ($(\angb:\raycone)+(\decal,0)$) -- (\decal,0);
\draw[dotted,thin,\colorboundcone] ($(\angc:\raycone)+(0,\decal)$) -- (0,\decal);

\end{tikzpicture}
\caption{$\cR_x$}

\end{subfigure}
\begin{subfigure}{.2\textwidth}
    
\begin{tikzpicture}[scale=\scalecRx]

\pgfmathsetmacro\anga{0};
\pgfmathsetmacro\angb{45};
\pgfmathsetmacro\angc{90};
\pgfmathsetmacro\angd{150};
\pgfmathsetmacro\ange{170};

\pgfmathparse{\anga>\angb}%
        \ifnum\pgfmathresult=1
                \pgfmathsetmacro\anga{\anga+360}%
        \fi
\fill[\colorcone,opacity=\opacitycone] ($(\anga:\raycone)+(\decal,0)$) -- (\decal,0) --  ($(\angb:\raycone)+(\decal,0)$)  arc (\angb:\anga:\raycone);

\pgfmathparse{\angb>\angc}%
        \ifnum\pgfmathresult=1
                \pgfmathsetmacro\angb{\angb+360}%
        \fi
\fill[\colorcone,opacity=\opacitycone] ($(\angb:\raycone)+(0,\decal)$) -- (0,\decal) -- ($(\angc:\raycone)+(0,\decal)$)  arc (\angc:\angb:\raycone);

\draw[thick,->,\colorboundcone,name path=halflinea] (\decal, 0) -- ($(\anga:\lengthcone)+(\decal,0)$);
\draw[thick,->,\colorboundcone,name path=halflineb] (\decal,0) -- ($(\angb:\lengthcone)+(\decal,0)$);

\draw[thick,->,\colorboundcone,name path=halflineb] (0,\decal) -- ($(\angb:\lengthcone)+(0,\decal)$);
\draw[thick,->,\colorboundcone,name path=halflinec] (0,\decal) -- ($(\angc:\lengthcone)+(0,\decal)$);

\draw[\colorvrtx] (\decal,0)node{$\bullet$};
\draw[\colorvrtx] (0,\decal)node{$\bullet$};

\end{tikzpicture}
\caption{$\bar \cR_x$}

\end{subfigure}
\begin{subfigure}{.2\textwidth}
    
\begin{tikzpicture}[scale=\scalecRx]

\pgfmathsetmacro\anga{0};
\pgfmathsetmacro\angb{45};
\pgfmathsetmacro\angc{90};
\pgfmathsetmacro\angd{150};
\pgfmathsetmacro\ange{170};

\pgfmathparse{\anga>\angb}%
        \ifnum\pgfmathresult=1
                \pgfmathsetmacro\anga{\anga+360}%
        \fi
\fill[\colorcone,opacity=\opacitycone] ($(\anga:\raycone)+(\decal,0)$) -- (\decal,0) --  ($(\angb:\raycone)+(\decal,0)$)  arc (\angb:\anga:\raycone);

\pgfmathparse{\angb>\angc}%
        \ifnum\pgfmathresult=1
                \pgfmathsetmacro\angb{\angb+360}%
        \fi
\fill[\colorcone,opacity=\opacitycone] ($(\angb:\raycone)+(0,\decal)$) -- (0,\decal) -- ($(\angc:\raycone)+(0,\decal)$)  arc (\angc:\angb:\raycone);

\draw[thick,->,\colorboundcone,name path=halflinea] (\decal, 0) -- ($(\anga:\lengthcone)+(\decal,0)$);

\draw[thick,->,\colorboundcone,name path=halflineb] (0,\decal) -- ($(\angb:\lengthcone)+(0,\decal)$);
\draw[thick,->,\colorboundcone,name path=halflinec] (0,\decal) -- ($(\angc:\lengthcone)+(0,\decal)$);

\draw[\colorvrtx] (0,\decal)node{$\bullet$};

\draw[dotted,thin,\colorboundcone] ($(\angb:\raycone)+(\decal,0)$) -- (\decal,0);

\end{tikzpicture}
\caption{$\cP$}

\end{subfigure}
\begin{subfigure}{.2\textwidth}
    
\begin{tikzpicture}[scale=\scalecRx]

\pgfmathsetmacro\anga{0};
\pgfmathsetmacro\angb{45};
\pgfmathsetmacro\angc{90};
\pgfmathsetmacro\angd{150};
\pgfmathsetmacro\ange{170};

\pgfmathparse{\anga>\angb}%
        \ifnum\pgfmathresult=1
                \pgfmathsetmacro\anga{\anga+360}%
        \fi
\fill[\colorcone,opacity=\opacitycone] ($(\anga:\raycone)+(\decal,0)$) -- (\decal,0) --  ($(\angb:\raycone)+(\decal,0)$)  arc (\angb:\anga:\raycone);

\pgfmathparse{\angb>\angc}%
        \ifnum\pgfmathresult=1
                \pgfmathsetmacro\angb{\angb+360}%
        \fi
\fill[\colorcone,opacity=\opacitycone] ($(\angb:\raycone)+(0,\decal)$) -- (0,\decal) -- ($(\angc:\raycone)+(0,\decal)$)  arc (\angc:\angb:\raycone);

\draw[thick,->,\colorboundcone,name path=halflinea] (\decal, 0) -- ($(\anga:\lengthcone)+(\decal,0)$);
\draw[thick,->,\colorboundcone,name path=halflineb] (\decal,0) -- ($(\angb:\lengthcone)+(\decal,0)$);

\draw[thick,->,\colorboundcone,name path=halflinec] (0,\decal) -- ($(\angc:\lengthcone)+(0,\decal)$);

\draw[\colorvrtx] (\decal,0)node{$\bullet$};

\draw[dotted,thin,\colorboundcone] ($(\angb:\raycone)+(0,\decal)$) -- (0,\decal);

\end{tikzpicture}
\caption{$\cP'$}

\end{subfigure}
\centering
\caption{$\cR_x$ is a partition of $\Xi$ in $6$ elements, $\bar \cR_x$ is not a partition, $\cP$  and $\cP'$ are two distinct coarsest partitions (in $2$ elements) with $\cR_x\refleq \cP \refleq \bar \cR_x$ and $\cR_x\refleq \cP' \refleq \bar \cR_x$.
\label{fig:Ra_and_Rbarx}}
\end{figure}
\begin{remark}
	When the distribution of $\va \xi$ is absolutely continuous 
	with respect to the Lebesgue measure of $\Xi$, $\cR_x \sim_\PP \bar \cR_x$, thus  $\cR_x$ is the coarsest partition adapted to $x$.

	If $\va \xi$ does not admit a density, 
	$\cR_x$ is still an adapted partition but not necessarily the coarsest, 
	which might not exist (see \cref{fig:Ra_and_Rbarx}).
	Nevertheless, any adapted partition should refine $\bar{\cR}_x$. 
	Unfortunately, we cannot use $\bar{\cR}_x$ in \GGAPM, 
	as we cannot guarantee that $\bar{\cR}_x$ is a $\PP$-partition.
\end{remark}

\begin{remark}
	Note that \cref{thm:ref_equiv_exact} implies that all user-defined partitions 
	$\cP^{(t)}$ of algorithm 1. of \cite{ramirez2021generalized} must satisfy $\cP^{(t)} \refleq_\PP \bar \cR_x$ by their proposition 2.
	In the finite scenario case, our adaptedness condition 
	is equivalent to Song and Luedtke's condition \cite{song2015adaptive}.
\end{remark}

We start the proof with a technical lemma. 
\begin{lemma}
	Consider a set $P \subset \Xi$ such that $\PP(P)>0$, and a state $x \in \RR^n$. 
	Then,
	\begin{align*}
	\exists R \in & \cR_x, \; P \subset_\PP R  \\
	&\Longrightarrow Q(x,\besp{\va \xi|P})=\besp{Q(x,\va \xi)|P} \\
	\exists \bar R \in &  \bar \cR_x,  \; P \subset_\PP \bar R  \\ 
	&\iff Q(x,\besp{\va \xi|P})=\besp{Q(x,\va \xi)|P}
	\end{align*}
	\label{lem:ref_Q_equality}
\end{lemma}

\begin{proof}
	Since $\exists R \in  \cR_x, \; P \subset_\PP R $ implies $\exists \bar R \in  \bar \cR_x, \; P \subset_\PP \bar R$, we only need to prove the second equivalence.

	($\Rightarrow$)
	Let $P$ be such that there exists $N \in \cN(D)$ with $P \subset_\PP \bar E_{N,x}$. 
	As $N\in \cN(D)$, there exists a dual point $\dualvar_N \in D$ such that $N$ is the normal cone of $D$ at $\lambda_N$.
	By definition of a normal cone, for all $\tilde h \in N$ and all $\lambda \in D$, $\tilde h^\top (\lambda - \lambda_N) \leq 0$. 
	In other words $\tilde h^\top \lambda_N =\max_{\lambda \in D} \tilde h^\top \lambda$.

	As $P\subset \bar E_{N,x}$, for $\PP$-almost-all $\xi \in P$, we have $h - Tx \in N$.
	Recall that  $Q(x,\xi)= \sup_{\dualvar \in D} (h-Tx)^\top \dualvar$, 
	thus, 
	 $Q(x,\xi)= (h-T x)^\top \dualvar_N$.
	 Hence,
	\begin{subequations}
	\begin{align*}
		\besp{Q(x,\va \xi)|P}
		&=\besp{(\va h-\va T x)^\top \dualvar_N|P} \\
		&=\besp{\va h-\va T x|P}^\top\dualvar_N \\
		&=Q(x,\besp{\va \xi|P})
	\end{align*}
	\end{subequations}
	 as $N$ is convex and $\besp{\va h -\va Tx |P} \in N$.

	($\Leftarrow$)
	For $\psi \in \RR^l$, we denote 
	$D^\psi:= \argmax_{\lambda \in D} \psi^\top \lambda$.
	Note that, for all $\psi, \psi' \in \relint(N)$, with $N \in \cN(D)$, we have $D^N:=D^\psi = D^{\psi'}$.

		Assume that there is no $R \in \bar \cR_x$ such that $P \subset_\PP R$.
Then, for all $R \in \bar \cR_x $, $\bprob{P \cap R}<\bprob{P}$.
	Since $\bprob{P}\leq \sum_{R \in \bar \cR_x} \bprob{P \cap R}$ 
	there exist $R_1$ and $R_2$ in $\bar \cR_x$ such that 
	$\bprob{P \cap R_1}>0$ and $\bprob{P \cap R_2}>0$.
		Let $\dualvar \in D$ such that $Q(x,\besp{\va \xi|P})=\besp{\va h -\va T x|P}^\top \dualvar$ i.e. $\dualvar \in D^{\besp{\va h -\va T x|P}}$.
		Let $N_1$ and $N_2 \in \collecmax{\cN(D)}$ be such that $R_1=\bar E_{N_1,x}$ and $R_2=\bar E_{N_2,x}$.
		Since $N_1 \neq N_2$ are maximal, $D^{N_1} \cap D^{N_2}=\emptyset$.
		Thus, there exists at least one $i \in \{1,2\}$ such that $\dualvar \not \in D^{N_i}$.
		Then, $\besp{Q(x,\va \xi)|P\cap R_i}>\besp{\va h - \va T x|P\cap R_i}^\top \lambda$.

		Note that $Q(x,\xi)=\sigma_{D}(h-Tx)$, 
		where $\sigma_D$ is the support function of the polyhedron $D$,
		thus $\xi \mapsto Q(x,\xi)$ is a polyhedral function.
		Further, its affine regions are the elements of $\bar \cR_x$.
		
		By convexity, for any measurable set $A$, 
		$\besp{Q(x,\va \xi)| P \cap A } \geq Q\bp{x,\besp{\va \xi| P \cap A }}$ which is equal to $\max_{\dualvar' \in D}\besp{\va h -\va Tx| P\cap A }^\top \dualvar'$. 
		Since $\dualvar \in D$, we have $\besp{Q(x,\va \xi)| P \cap A }\geq \besp{\va h -\va Tx| P\cap A }^\top \dualvar$.

		Thus, $\besp{Q(x,\va \xi)| P } > Q\bp{x,\besp{\va \xi| P }} $.

\end{proof}

\begin{proof}[Proof of \cref{thm:ref_equiv_exact}]

By definition $\cP \refleq_\PP \bar \cR_x$
if and only if for all $P\in \cP$ there exists a cell $R \in \bar \cR_x$ such that $ P \subset_\PP R$.
By~\cref{lem:ref_Q_equality} this is equivalent to, for all $P\in \cP$, 
$Q(x,\besp{\va \xi|P})= \besp{Q(x,\va \xi)|P}$.
Now by Jensen's inequality, this equality for all $P\in \cP$ is equivalent to the equality of a convex sum like
$$ \sum_{P \in \cP} Q\bp{x,\besp{\va \xi|P}} \bprob{P} =\sum_{P \in \cP} \besp{Q(x,\va \xi)|P} \bprob{P} .$$
Law of total expectation yields \eqref{eq:ref_equiv_exact}.

\end{proof}

\begin{remark}
	Let $x^\star$ be an optimal solution of
	\begin{equation}
			\min_{x \in X} \quad c^\top x  +\ectg_{\cP^\star}(x) 
	\end{equation}
	where $\cP^\star \preceq_\PP \cR_{x^\star}$.
	Then, $x^\star$ is also a solution of Problem \eqref{eq:2SLP}.
\end{remark}

\subsection{Generalization with stochastic second cost q}

Up to now, $\va h$ and $\va T$ have arbitrary distribution.
In this subsection, we explain how,
mainly from a theoretical point of view,
 \GGAPMspace can be extended to 2SLP, with arbitrary cost distribution and finitely supported $\va W$.

\subsubsection{Finitely supported W and q}

When the cost $\va q$ and the recourse matrix $\va W$ are not fixed, the dual admissible set $D$ depends on the scenario.
We denote
$D_{W,q}:=\{\dualvar \in \RR^l \, |\, W^\top \dualvar \leq q\}$.
We cannot work with a constant normal fan $\cN(D)$ and apply directly \GGAPM.
Nevertheless, 
	by defining $\tilde \ectg(x|W,q):=\besp{Q(x,\va \xi)|(\va W,\va q)=(W,q)}$,
	we can adapt \GGAPMspace by finding an adapted partition $\cR_{x,W,q}$ for each $\tilde \ectg(.|W,q)$ and taking the common refinement $\cR_x=\bigwedge_{(W,q) \in \supp(\va W, \va q)}\cR_{x,W,q}$
	Then, we solve the master problem with $\ectg_{\cP}(x):=\sum_{(W,q) \in \supp(\va W,\va q)}\bprob{W,q} \tilde \ectg_{\cP}(x|W,q)$.

\subsubsection{General stochastic second cost q}

We now consider a cost $\va q$ with arbitrary distribution.
Denote $\secfan{W}$ the secondary fan of $W$ 
(see e.g. \cite[Chapter 5]{de2010triangulations}), and $\cQ_W = \ba{ \relint(S) \; | \; S \in  \secfan{W} }$. 
Note that $\cQ_W$ is a finite partition of  the space of $q$ such that 
 $q \mapsto \cN(D_{W,q})$ is constant in each cell of $\cQ_W$ (see \S~3.2 in \cite{forcier2021exact}), and so is $\cR_{x,W,q}$. 

 	Thus, for a given $W$, the common refinement $\bigwedge_q \cR_{x,W,q})$ is actually the common refinement of a finite number of partitions.
 	Finally, as $\supp(\va W)$ is finite,
	$\cR_x = \bigwedge_{(W,q) \in \supp(\va W, \va q)} \cR_{x,W,q}$ is still a well-defined and finite collection for every distribution of $q$.
	We can still solve the master problem with
	$\ectg_{\cP}(x):=\sum_{W\in \supp(\va W)}\sum_{R\in \cQ_W}\bprob{\va W=W,\va q \in R} \tilde \ectg_{\cP}(x|\va W=W,\va q \in R)$.

However, we cannot apply this method to non finitely supported $\va W$.

\section{Comparison with other algorithms and convergence}
\label{sec:comparison_convergence}

In this section, we show that the partition based methods  can be seen as an acceleration of the cutting plane method.
It then gives us a finite convergence proof with a bound on the number of steps.

\subsection{Adapted partition and subdifferential}

We show that, for any first stage control $x \in X$,
if the partition is adapted to $x$, 
then the subdifferential of aggregated recourse cost coincide with the subdifferential of the true recourse cost.

\begin{lemma}
	Let $x \in \dom(\ectg)$ and $\cP$ be a refinement of $\cR_x$, i.e. $\cP\refleq \cR_x$, then
	\begin{equation}
		\partial \ectg_{\cR_x}(x) \subset \partial \ectg_\cP(x) \subset \partial \ectg(x)
	\end{equation}
	Furthermore, if $x \in \relint \dom(\ectg)$,
	\begin{equation}
		\partial \ectg_{\cR_x}(x) = \partial \ectg_\cP(x) =\partial \ectg(x)
	\end{equation}
	\label{lem:subgradient_ectg_cP}
\end{lemma}

\begin{proof}
  
	Let $g \in \partial \ectg_{\cR_x}(x)$ then for all $y$, $\ectg_{\cR_x}(y) \geq \ectg_{\cR_x}(x) +g^\top (y-x)$, by monotonicity $\ectg_\cP(y) \geq \ectg_{\cR_x}(y)$ 
	and by construction of $\cR_x$, $\ectg_{\cR_x}(x)=\ectg(x)=\ectg_\cP(x)$. Thus, $\ectg_{\cP}(y) \geq \ectg_{\cP}(x)+ g^\top (y-x)$ and $g \in \partial \ectg_\cP(x)$. The proof for the second inclusion is similar.

	Let $x \in \relint \dom (\ectg)$, we now prove that $\partial \ectg_{\cR_x}(x)=\partial \ectg(x)$.
	Recall that $D^N=D^\psi= \argmax_{\lambda \in D} \psi^\top \lambda$, for $\psi \in \relint(N)$ where $N\in \cN(D)$.
	By \cite[Thm. 11 p. 117]{birge2011introduction}, 
	$\partial \ectg(x)= \besp{\partial_x Q(x,\va \xi)} + N_{\dom(V)}(x)$ and by \cite[Prop 2.2 p.28]{shapiro2014lectures}, $\partial_x Q(x,\xi)=-T^\top D^{h-Tx}$. Thus, since $x \in \relint \dom (\ectg)$,
	\begin{subequations}
		\begin{align}
			\partial \ectg(x)
			& = \besp{-\va T^\top D^{\va h - \va T x}} \\
			&= \besp{\sum_{N \in \cN(D)} -\indi{\va h-\va Tx \in \relint(N)}\va T^\top D^N } \\
			&= \besp{\sum_{N \in \cN(D)} -\indi{\va \xi \in E_{N,x}}\va T^\top D^N }
		\end{align}
	
Further, 
	\begin{align*}
		\EE \big[ & \indi{\va \xi \in E_{N,x}}  \va T^\top D^N  \big]\\ 
		& = \bprob{E_{N,x}}\besp{\va T|E_{N,x}}^\top D^N\\
& = \bprob{E_{N,x}}\besp{\va T|E_{N,x}}^\top D^{\besp{\va h -\va T x|E_{N,x}}}  \\
	\end{align*}
And by definition of $\cR_x$ we get 
	\begin{align}
		\partial \ectg(x) &= \sum_{P \in \cR_x} -\bprob{P}\besp{\va T|P}^\top D^{\besp{\va h -\va T x|P}}  \\
			&= \sum_{P \in \cR_x} \bprob{P}\partial_x Q(x,\besp{\va \xi|P})  \\
			&= \partial \ectg_{\cR_x}(x)
	\end{align}
	\end{subequations}
\end{proof}

\subsection{Link with L-shaped and Bender decomposition}

The classical L-shaped method (see e.g. \cite{birge2011introduction}) is a specification of Bender's decomposition to 2SLP for finitely supported distribution. 
The core idea consists in representing the expected cost-to-go function in~\eqref{eq:2SLP}, by a lift variable
\begin{equation}
		\min_{x\in X, \theta \in \RR} \; \Ba{c^\top x + \theta \; | \; (x,\theta) \in \epi (\ectg)} 
\label{eq:2SLP_lift}
\end{equation}
We then relax the epigraphical representation 
$(x,\theta) \in \epi (\ectg)$, replacing it by a set of valid inequalities called \emph{cuts}, i.e.
\begin{subequations}
	\begin{align}
		\min_{x\in X, \theta \in \RR} \quad &  c^\top x + \theta \\
		\st \quad & 
		g^\top x + v \leq \theta, & \forall (g,v)\in \cO, \\
		 & f^\top x \leq \bar{f}, & \forall (f,\bar{f}) \in \cF .
	\end{align}
	\label{eq:master_L_shaped}
\end{subequations}

More precisely, assume that we have such a relaxation of~\eqref{eq:2SLP}.
Let $x^k$ be an optimal first stage control of this relaxation.
If it is admissible, meaning that for all scenario $\xi$ there exists an admissible recourse control $y_\xi$, we compute, through duality, a subgradient $g^k \in \partial \ectg(x^k)$. 
This yields
a new \emph{optimality cut} $\theta \geq (g^k)^\top(x - x^k) + \ectg(x^k)$, which is added to $\cO$.
If $x^k$ is not admissible we can add a \emph{feasibility cut} to $\cF$ instead.
We then solve our strengthened relaxation to obtain $x^{k+1}$. 

The L-Shaped method specifies that the subgradient $g^{k}$ can be obtained as an average over $\xi$ of subgradients $g^{k,\xi} \in \partial_x Q(x^{k},\xi)$.
In particular, it means that, to compute the subgradient, we can solve $l$ smaller LP instead of a large one.

\begin{remark}[L-shaped for continuous distribution]
	When the distribution are non-finitely supported, we cannot apply naively this method as there is a non-finite number of scenarios.
Nevertheless,  we can still approximate $\epi(\ectg)$ with cuts.
We can compute $\theta=\ectg_{\cR_x}(x)$ and a subgradient $g \in \partial \ectg_{\cR_x}(x)$ by solving $|\cR_x|$ linear problems of the form \eqref{eq:dual_2nd_stage} through exact quantization.
By \cref{thm:ref_equiv_exact}, $\theta=\ectg_{\cR_x}(x)=\ectg(x)$.
Further,  $g \in \partial \ectg_{\cR_x}(x) \subset \partial \ectg(x)$ by \cref{lem:subgradient_ectg_cP}.
Then $(\theta,g)$ defined an optimality cut.
\end{remark}

\cref{lem:subgradient_ectg_cP} shows that, at each step $k$ of \GGAPMspace algorithm, we add a collection of valid cuts which are exact at $x_k$ to our collection of cuts.
This means that APM methods can be seen as a Bender's decomposition method where we add more than one exact cut per iteration.
In particular, when $x^k \in \relint \dom (\ectg)$ we add the whole tangent cone of $\epi(\ectg)$ at $x$ instead of a single cut.

\subsection{Convergence of \GGAPM}
We start by showing that the bounds generated in \cref{algo:GAPM} are monotonous.
\begin{lemma}
\label{lem:bounds_monotonicity}
	For every computed step $k$ we have
	\begin{equation}
		z_L^{k-1} \leq z_L^{k} \leq \val\eqref{eq:2SLP} \leq z_U^{k} \leq z_U^{k-1}
	\end{equation}
\end{lemma}

\begin{proof}
	 Since $\cP^{k}\refleq_\PP \cP^{k-1}$, by \cref{lem:monotony_ref}, we have, for all $x \in X$,
	\begin{subequations}
		\begin{align}
		c^\top x + \ectg_{\cP^{k-1}}(x)
		&\leq c^\top x +  \ectg_{\cP^{k}}(x) \\
		&\leq  c^\top x + \ectg(x)
		\end{align}
	\end{subequations}
	Minimizing over $x$ yields $z^L_{k-1} \leq z^L_{k} \leq \val\eqref{eq:2SLP}$.
	$ z^U_{k} \leq z^U_{k-1}$ is direct from the definition.

	Note that $ z^U_{0} = +\infty \geq \val\eqref{eq:2SLP}$.
	Assume that $z^U_{k-1} \geq \val\eqref{eq:2SLP}$. 
	If $z^U_{k}<z^U_{k-1}$, 
	as $\cP^{k} \refleq_\PP \cR_{x^{k}}$ we have
	\begin{equation}
		z^U_{k}=c^\top x_{k} + \ectg_{\cP^k}(x_k)=  c^\top x_{k}+ \ectg(x_k)
	\end{equation}
	which is an upper bound of $\val\eqref{eq:2SLP}$ since $x_{k} \in X$.
	Induction yields $z^U_{k} \geq \val\eqref{eq:2SLP}$ for all $k$.
\end{proof}

We can now prove finite convergence of any APM method.

\begin{theorem}
	If $X \cap \dom(\ectg) \subset \RR^+$ has a finite diameter $M \in \RR^+$ and $x\to c^\top x+\ectg(x)$ is Lipschitz with constant $L$ then the partition based \cref{algo:GAPM} finds an $\varepsilon$-solution in at most $\bp{\frac{\sqrt n LM}{\varepsilon}+1}^n$ iterations.
\end{theorem}

\begin{proof}
	We adapt the classical proof of Kelley's cutting plane algorithm to APMs.
	Let $k \in \NN$ and $1<i<k$, we have that $ \ectg(x_i)= \ectg_{\cP^{k-1}}(x_i)= \ectg_{\cP^{i}}(x_i)$.
	Let 
	 $g \in \partial \ectg_{\cP^{k-1}}(x_i) \subset \partial \ectg(x_i)$ such that $\|c+g\|$ is bounded by the Lipschitz constant $L$ then
	\begin{subequations}
		\begin{align*}
			z^U_k-z_k^L
			&\leq c^\top x_i+\ectg_{\cP^i}(x_i)-\bp{c^\top x_k+\ectg_{\cP^{k-1}}(x_k)} \\
			&=c^\top (x_i-x_k)+\ectg_{\cP^{k-1}}(x_i)-\ectg_{\cP^{k-1}}(x_k)) \\
			&\leq c^\top (x_i-x_k) -g^\top(x_k-x_i) \\
			&\leq \|c+g\|_2 \|x_i-x_k\|_2 \leq L\|x_i-x_k\|_2
		\end{align*}
	\end{subequations}
	Then, for $k$ such that, $\varepsilon < z^U_k-z_k^L$, we have $\varepsilon <\sqrt n L\|x_i-x_k\|_\infty$, in particular $\|x_i-x_k\|_\infty \geq \varepsilon/(\sqrt n L)$.
	By definition of $M$ there are at most
	$\bp{\frac{\sqrt n LM}{\varepsilon}+1}^n$ boxes of radius $ \varepsilon/(\sqrt n L)$.
	An $\varepsilon$-solution being obtained as soon as two points are in the same box.
\end{proof}

\section{Numerical examples}
\label{sec:numerical}

In this section, we detail the actual computation required by \GGAPMspace and illustrate the algorithm on numerical examples.

\subsection{Detailing computation}

In the following two sections, we give more details on how to compute the \cref{li:solve_part_lp,li:adapted_partition,li:common_refinement,li:upper_bound} of  \cref{algo:GAPM}.

\subsubsection{Master problem and subproblems}

Once $\besp{\va \xi \, |\, P}$ and $\bprob{P}$ have been are computed for $P \in \cP^{k-1}$, by \cref{eq:defi_ectg_partition} and \cref{eq:2nd_stage}, the problem  of \cref{li:solve_part_lp} is reduced to the following linear problem

\begin{subequations}
	\begin{align*}
		\min_{x \in X, (y_P) \in (\RR_+^m)^{\cP^{k-1}}} &  c^\top x + \sum_{P \in \cP^k} \bprob{P} q^\top y_P \\
		\st \qquad & \besp{\va T|P} x + W y_P = \besp{\va h|P} \nonumber  \\
		& \hspace{2cm} \forall P \in \cP^k. 
	\end{align*}
\end{subequations}

Moreover, to compute the upper bound in \cref{li:upper_bound}, we need to solve at most $|\cP^k|$ linear problems of dimension $m$ 
\begin{subequations}
	\begin{align*}
		Q(x,\besp{\va\xi|P}):= \min_{y_P \in \RR_+^m} & \; q^\top y_P \\
		\st &\; \besp{\va T|P} x + W y_P = \besp{\va h|P} 
	\end{align*}
\end{subequations}

\subsubsection{Refinement, expectation and probabilities}

\label{sec:compute_besp_prob_ref}

Recall that we can store a polyhedron $E$, either as a family of constraints $(A,b)$ such that 
$E=\{x\in \RR^d\,|\, Ax\leq b\}$ ($H$-representation) or as families of vertices $(v_i)_{i\in I}$ and rays $(r_j)_{j \in J}$ such that $E=\conv(v_i)_{i\in I}+\cone(r_j)_{j\in J}$ ($V$-representation). 
Both representation are implemented in polymake \cite{polymake:2000}.
We can switch between representations through algorithms such as the \emph{double description} \cite{FukudaProdon96}.

We can simultaneously compute conditional expectation, probabilities and refinement as detailed in \cref{algo:common_ref_besp_prob}.

\begin{algorithm}
	\KwData{$\cP^{k-1}$ and $\cR_{x_k}$ the partition to refine, second stage distributions $\va T$ and $\va h$.}
	Set $\cP^{k}:=\emptyset$\;
		\For{$P \in \cP^{k-1}$ and $R \in \cR_{x^{k}}$}
		{
			Set $P':=P \cap R$\;
			\If{$\bprob{P'}>0$}
			{
				Store $\bprob{P'}$, $\besp{\va T | P'}$ and $\besp{\va h | P'}$\; \label{li:store_prob_esp}
				Set $\cP^{k}:=\cP^{k} \cup \{P'\}$\;
			}
		}
	\caption{Computing refinement}
	\label{algo:common_ref_besp_prob}
\end{algorithm}

In this algorithm, the computation of probabilities on polyhedra in \cref{li:store_prob_esp} is a $\sharp$-P complete problem in the general case.
For a large class of distributions, formulas exists. 
See \cite[Section 5]{forcier2021exact} for a review.

\subsubsection{Computing the adapted partition}

\label{sec:compute_adapted_partition}

In this section, we explain how to compute $\cR_x = \ba{ E_{N,x} \, |\, N \in \cN(D)}$ where $ E_{N,x} = \ba{(T,h) \,|\, h-Tx \in \relint(N)}$.

The normal fan $\cN(D)$ computation, already implemented in polymake, can be done thanks to a double description and active constraint sets.
Note that if $N \in \cN(D)$, then $E_{N,x}$ is a relatively open polyhedral cone of $\Xi$.
In particular, if $N:= \{\tilde h \, |\, M \tilde h \leq 0\}$ is given in a non-redundant $H$-representation where $M \in \RR^{p \times l}$, we have $\relint(N)=\{\tilde h \, |\, M \tilde h \vecll 0\}$.
Then $E_{N,x}=\{\xi\in \Xi \, |\, H^x\xi \vecll 0\}$, with  $H^x=(- x_1 M \cdots -x_n M \; M)$.

Unfortunately, as normal cone are given in $V$-representation, the above approach require a double-description to obtain the $H$-representation of $N$.
This work well in low dimension, but is numerically intractable in higher dimension.
Indeed, McMullen \cite{mcmullen1970maximum} showed that the size of the $H$-representation can be exponential in the size of the $V$-representation and vice versa.

The double-description can be avoided if the technology matrix $\va T\equiv T$ is fixed.
Indeed, in this case
$ E_{N,x} \sim_\PP \{T\} \times \bp{Tx + \relint(N)}$.

Thus, we can compute, at the beginning of the algorithm, 
a  $V$-representation of all $N \in \cN(D)$,
and easily deduce a $V$-representation of $E_{N,x}$ by adding $Tx$ to each representant ray.

\subsection{Numerical examples}

The problems presented in \cite{ramirez2021generalized} either have a fixed technology matrix or a random variable with a support of dimension 1.
We applied \GGAPM~to the problems LanDs and CV@R (ibid.) and to the problem Prod-Mix that cannot be tackled directly by GAPM. Our code is available at
	\url{https://github.com/maelforcier/GAPM}

\subsubsection{Energy planing problem - LandS}
We applied numerically our methods to the LandS problem and found the same partition, state $x_k$, lower and upper bounds as in \cite{ramirez2021generalized}.

\subsubsection{Conditional value-at-risk linear problems}

For the conditional value-at-risk problem in \cite{ramirez2021generalized}, note that our analysis through normal fans leads to the same partition:
\begin{subequations}
	\begin{align}
		Q^D(x,\xi):= \quad \max_{\lambda \in \RR}
		& \quad (-x^\top r^\xi -\tau)\lambda \\
		\st & \quad 0\leq \lambda \leq 1
	\end{align}
\end{subequations}
Here $D=[0,1]$ and $\cN(D)=\{\RR^-,\{0\},\RR^+\}$
Then, if $x\neq 0$,
$\cR_x = \ba{\{r |x^\top r>-\tau \}, \{r |x^\top r=-\tau \}, \{r |x^\top r<-\tau \}}$.

\subsubsection{Prod-Mix}
\label{sssec:prod_mix} 
We adapted the problem Prod-mix of \url{https://stoprog.org/SavedLinks/IBM_StoExt_problems/node4.php} as
\begin{align*}
	\min_{x,\va y} &\quad -c^\top x+\besp{q^\top \va y} \\
	\st &\quad  \va Tx-\va y \leq \va h \\
	& \quad x,\va y\geq 0,
\end{align*}
where $q^\top=(5,10)$, $c^\top=(12,40)$, $\va T$ follows the uniform law $\begin{pmatrix}
	\cU[3.5,4.5] & \cU[9,11] \\
	\cU[0.8,1.2] & \cU[36,44]
\end{pmatrix}$
and $\va h^\top$ follows the uniform distribution $\bp{\cU[5970,6030], \;\cU[3979,4021]}$. 
\GGAPMspace gave the results summed up in \cref{tab:result_prod_mix}




\begin{table}[h]
\centering
\footnotesize{
\begin{tabular}{|c|c|c|c|c|}
\hline
	 $k$ &$x_k$ & $z_L^k$ & $z_U^k$ & $|\collecmax{\cP_k}|$\\
	 \hline
	$1$ & $(1333.33, 66.67)$ & $-18666.67$ & $-16939.71$ & $4$\\\hline
	$4$ & $(1379.98, 56.64)$ & $-17744.67$ & $-17708.00$ & $25$\\\hline
	$7$ & $(1379.31, 55.82)$ & $-17712.06$ & $-17711.45$ & $64$\\\hline
	$10$& $(1377.69, 56.05)$ & $-17711.57$ & $-17711.56$ & $121$\\\hline
\end{tabular}
}
\caption{Results of \GGAPMspace for Prod-Mix \label{tab:result_prod_mix}}
\end{table}

To compare our approach with SAA, we solved the same problem 100 times, each with 10 000 scenarios randomly drawn, yielding a $95\%$ confidence interval centered in $-17711$, with radius $2.2$.
Thus, \GGAPMspace can be useful to find accurate values.

The most time-consuming parts of the algorithm
are the computations of volumes, because polymake only implement exact computations, which was proven to be $\sharp P$-complete \cite{dyer1988complexity}.
To improve \GGAPM, we could use precise rapid approximation volume algorithms, see e.g. \cite{cousins2016practical}.

\bibliographystyle{alpha}

\end{document}

%% file: commandes.tex
\newcommand{\ectg}{V} 

\newcommand{\secfan}[1]{\Sigma\operatorname{-fan}(#1)}

\newcommand{\dualvar}{\lambda}

\newcommand{\EE}{{\mathbb E}} 

\newcommand{\NN}{{\mathbb N}} 
\newcommand{\PP}{{\mathbb P}} 
\newcommand{\RR}{{\mathbb R}} 

\newcommand{\cA}{\mathcal{A}}

\newcommand{\cF}{\mathcal{F}}

\newcommand{\cN}{\mathcal{N}}
\newcommand{\cO}{\mathcal{O}}
\newcommand{\cP}{\mathcal{P}}
\newcommand{\cQ}{\mathcal{Q}}
\newcommand{\cR}{\mathcal{R}}

\newcommand{\cU}{\mathcal{U}}

\newcommand{\indi}[1]{\mathds{1}_{#1}}

\newcommand{\collecmax}[1]{#1^{\max}}
\newcommand{\cone}{\operatorname{Cone}}
\newcommand{\conv}{\operatorname{Conv}}
\newcommand{\dom}{\operatorname{dom}}
\newcommand{\epi}{\operatorname{epi}}
\newcommand{\ie}{\emph{i.e.}}

\newcommand{\relint}{\operatorname{ri}}
\newcommand{\st}{\text{s.t. }}
\newcommand{\supp}{\operatorname{supp}}

\newcommand{\val}{\operatorname{val}}

\newcommand{\refleq}{\preceq}

\newcommand{\vecll}{\ll}

\renewcommand{\geq}{\geqslant}
\renewcommand{\leq}{\leqslant}
\renewcommand{\preceq}{\preccurlyeq}

{\normalsize}

{\normalsize}

{\normalsize}

\renewcommand{\bar}{\overline}
\renewcommand{\tilde}{\widetilde}

\newcommand{\defegal}{:=}

\newcommand{\bp}[1]{\big(#1\big)}                           
\newcommand{\Bp}[1]{\Big(#1\Big)}                           

\newcommand{\bc}[1]{\big[#1\big]}                           

\newcommand{\ba}[1]{\big\{#1\big\}}                         
\newcommand{\Ba}[1]{\Big\{#1\Big\}}                         

\newcommand{\argmax}{\mathop{\mathrm{argmax}}}

\newcommand{\prbt}{\mathbb{P}}                              
\newcommand{\mathscr}{\EuScript}

\newcommand{\va}[1]{\bm{#1}}

\newcommand{\bprob}[1]{\prbt\bc{#1}}                  

\newcommand{\espe}{\mathbb{E}}                              

\newcommand{\besp}[2][]{\espe_{#1}\bc{#2}}            

%% file: main.bbl
\begin{thebibliography}{test}

\bibitem[Bir85]{birge1985decomposition}
John~R Birge.
\newblock Decomposition and partitioning methods for multistage stochastic
  linear programs.
\newblock {\em Operations Research}, 33(5):989--1007, 1985.

\bibitem[BL11]{birge2011introduction}
John~R Birge and Francois Louveaux.
\newblock {\em Introduction to stochastic programming}.
\newblock Springer Science \& Business Media, 2011.

\bibitem[CV16]{cousins2016practical}
Ben Cousins and Santosh Vempala.
\newblock A practical volume algorithm.
\newblock {\em Mathematical Programming Computation}, 8(2):133--160, 2016.

\bibitem[DF88]{dyer1988complexity}
Martin~E. Dyer and Alan~M. Frieze.
\newblock On the complexity of computing the volume of a polyhedron.
\newblock {\em SIAM Journal on Computing}, 17(5):967--974, 1988.

\bibitem[DGKR03]{dupacova2003scenario}
Jitka Dupa{\v{c}}ov{\'a}, Nicole Gr{\"o}we-Kuska, and Werner R{\"o}misch.
\newblock Scenario reduction in stochastic programming.
\newblock {\em Mathematical Programming}, 95(3):493--511, 2003.

\bibitem[DLRS10]{de2010triangulations}
Jes{\'u}s~A De~Loera, J{\"o}rg Rambau, and Francisco Santos.
\newblock {\em Triangulations Structures for algorithms and applications}.
\newblock Springer, 2010.

\bibitem[EZ94]{edirisinghe1994bounds}
NCP Edirisinghe and William~T Ziemba.
\newblock Bounds for two-stage stochastic programs with fixed recourse.
\newblock {\em Mathematics of Operations Research}, 19(2):292--313, 1994.

\bibitem[FGL21]{forcier2021exact}
Maël Forcier, Stéphane Gaubert, and Vincent Leclère.
\newblock Exact quantization of multistage stochastic linear problems, 2021.

\bibitem[FP96]{FukudaProdon96}
K.~Fukuda and A.~Prodon.
\newblock Double description method revisited.
\newblock In {\em Selected papers from the 8th Franco-Japanese and 4th
  Franco-Chinese Conference on Combinatorics and Computer Science}, pages
  91--111, London, UK, 1996. Springer-Verlag.

\bibitem[GJ00]{polymake:2000}
Ewgenij Gawrilow and Michael Joswig.
\newblock {\tt polymake}: a framework for analyzing convex polytopes.
\newblock In {\em Polytopes---combinatorics and computation ({O}berwolfach,
  1997)}, volume~29 of {\em DMV Sem.}, pages 43--73. Birkh\"auser, Basel, 2000.

\bibitem[HR07]{heitsch2007note}
Holger Heitsch and Werner R{\"o}misch.
\newblock A note on scenario reduction for two-stage stochastic programs.
\newblock {\em Operations Research Letters}, 35(6):731--738, 2007.

\bibitem[Kuh06]{kuhn2006generalized}
Daniel Kuhn.
\newblock {\em Generalized bounds for convex multistage stochastic programs},
  volume 548.
\newblock Springer Science \& Business Media, 2006.

\bibitem[McM70]{mcmullen1970maximum}
Peter McMullen.
\newblock The maximum numbers of faces of a convex polytope.
\newblock {\em Mathematika}, 17(2):179--184, 1970.

\bibitem[RPM21]{ramirez2021generalized}
Cristian Ramirez-Pico and Eduardo Moreno.
\newblock Generalized adaptive partition-based method for two-stage stochastic
  linear programs with fixed recourse.
\newblock {\em Mathematical Programming}, pages 1--20, 2021.

\bibitem[SDR14]{shapiro2014lectures}
Alexander Shapiro, Darinka Dentcheva, and Andrzej Ruszczy{\'n}ski.
\newblock {\em Lectures on stochastic programming: modeling and theory}.
\newblock SIAM, 2014.

\bibitem[SL15]{song2015adaptive}
Yongjia Song and James Luedtke.
\newblock An adaptive partition-based approach for solving two-stage stochastic
  programs with fixed recourse.
\newblock {\em SIAM Journal on Optimization}, 25(3):1344--1367, 2015.

\bibitem[vAdOS18]{van2018adaptive}
Wim van Ackooij, Welington de~Oliveira, and Yongjia Song.
\newblock Adaptive partition-based level decomposition methods for solving
  two-stage stochastic programs with fixed recourse.
\newblock {\em Informs Journal on Computing}, 30(1):57--70, 2018.

\bibitem[VSW69]{van1969shaped}
Richard~M Van~Slyke and Roger Wets.
\newblock L-shaped linear programs with applications to optimal control and
  stochastic programming.
\newblock {\em SIAM Journal on Applied Mathematics}, 17(4):638--663, 1969.

\bibitem[WZ05]{wallace2005applications}
Stein~W Wallace and William~T Ziemba.
\newblock {\em Applications of stochastic programming}.
\newblock SIAM, 2005.

\end{thebibliography}
